\documentclass[12pt,leqno]{amsart}
\usepackage{amssymb,amsthm}
\usepackage{amsmath,amssymb,color}
\newcommand{\ep}{\varepsilon}

\newcommand{\va}{\varphi}
\newcommand{\ppp}{\partial}

\newcommand{\R}{\mathbb{R}}

\newcommand{\N}{\mathbb{N}}

\newcommand{\DDD}{\mathcal{D}}
\newcommand{\www}{\widetilde}

\newcommand{\PPPP}{\mathcal{P}}
\newcommand{\AAAA}{\mathcal{A}}

\newcommand{\sumjk}{\sum_{j,k=1}^n}

\newcommand{\ooo}{\overline}
\newcommand{\OOO}{\Omega}

\title
[]
{
Inverse parabolic problems by Carleman estimates
with  data taken initial or final time moment of observation}

\author{
$^1$ O.Y. Imanuvilov and 
$^{2}$ M.~Yamamoto }

\thanks{
$^1$  
Department of Mathematics, Colorado State University, 101 Weber Building, 
Fort Collins CO 80523-1874, USA e-mail:{\tt oleg@math.colostate.edu}
\\
$^2$ Graduate School of Mathematical Sciences, The University
of Tokyo, Komaba, Meguro, Tokyo 153-8914, Japan 
e-mail: {\tt myama@ms.u-tokyo.ac.jp}
}

\date{}
\begin{document}
\maketitle

\baselineskip 18pt
\begin{center}
Dedicated to the memory of Professor Victor Isakov
\end{center}
\begin{abstract}
We consider a parabolic equation in a bounded domain $\OOO$
over a time interval $(0,T)$ 
with the homogeneous Neumann boundary condition.
We arbitrarily choose a subboundary $\Gamma \subset \ppp\OOO$.
Then, we discuss an inverse problem of determining a zeroth-order 
spatially varying coefficient by extra data of solution $u$:
$u\vert_{\Gamma \times (0,T)}$ and $u(\cdot,t_0)$ in $\OOO$
with $t_0=0$ or $t=T$.
First we establish a conditional Lipschitz stability estimate as well as
the uniqueness for the case $t_0=T.$
Second, under additional condition for 
$\Gamma$, we prove the uniqueness for the case $t_0=0$.
The second result adjusts the uniqueness by M.V. Klibanov
(Inverse Problems {\bf 8} (1992) 575-596) to the inverse problem 
in a bounded domain $\OOO$.  We modify his method which 
reduces the inverse parabolic problem to an inverse hyperbolic problem, and so 
even for the inverse parabolic problem,
we have to assume conditions for the uniqueness for the corresponding
inverse hyperbolic problem.
Moreover we prove the uniqueness for some inverse source problem for 
a parabolic equation for $t_0=0$ without boundary condition on the whole
$\ppp\OOO$.
\end{abstract}

%\begin{document}
%\maketitle

\section{Introduction and main results}

Let $\Omega$ be a bounded domain in $\Bbb R^n$ with $C^2$-boundary.
We fix a moment $t_0 \in [0,T]$, and introduce the elliptic operator
$$
Av(x) = \sum_{j,k=1}^n \partial_j(a_{jk}(x)\partial_k v(x))
- \sum_{j=1}^n b_j(x)\partial_jv(x) \quad \mbox{for $v\in H^2(\OOO)$}.
$$

We consider 
%\begin{equation}\label{Ssok1}
$$
\partial_tu = Au - p(x)u \quad 
\mbox{in $(0,T) \times \OOO$},          \eqno{(1.1)}
$$
$$
%\begin{equation}\label{Ssok2}
\partial_{\nu_A}u\vert_{(0,T) \times \ppp\OOO}=0,   \eqno{(1.2)}
$$
and
%\begin{equation}\label{Ssok3}
$$
u(t_0,\cdot)=u_0 \quad \mbox{in $\OOO$}.          \eqno{(1.3)}
$$

Here, choosing a constant $\kappa \in (0,1)$ arbitrarily, we assume 
%\begin{equation}\label{sok4}
$$
a_{jk}=a_{kj} \in C^{2+\kappa}(\ooo{\Omega}), \quad 
b_j, p\in C^{\kappa}(\ooo{\Omega}) \quad \mbox{for all 
$j,k\in \{1,\dots,n\}$}                 \eqno{(1.4)}
$$
and there exists a constant $\sigma>0$ such that
%\begin{equation}\label{sok5}
$$
\sum_{j,k=1}^n a_{jk}(x)\xi_j\xi_k\ge \sigma\vert \xi\vert^2\quad 
\mbox{for all $(x,\xi)\in \Omega\times \Bbb R^n$}.      \eqno{(1.5)}
$$
Here and henceforth, let $C^{2+\kappa}(\ooo{\OOO})$, $C^{\kappa}(\ooo{\OOO})$
with $\kappa \in (0,1)$ denote the Schauder spaces 
(e.g., Ladyzhenskaya, Solonnikov and Ural'ceva \cite{LSU}), and
let $\nu = (\nu_1, ..., \nu_n)$ be the unit outward
normal vector to $\ppp\OOO$.  We set
$$
\partial_{\nu_A} v = \sum_{k,j=1}^na_{kj}(x)\nu_j\partial_kv.   \eqno{(1.6)}
$$

Throughout this article, we fix the coefficients $a_{jk}, b_j$, 
$1\le j,k \le n$  in (1.1) and
by $u(c)=u(c)(x,t)$, we denote one solution to (1.1) - (1.3) with 
zeroth-order coefficient $c=c(x).$ 

In this article, we consider 
\\
{\bf Inverse coefficient problem.}
\\
{\it
Let $\Gamma$ be a fixed subboundary of $\partial\Omega$ and 
let $t_0 \in [0,T]$ be given in (1.3).
Determine $p(x)$ by $u\vert_{(0,T) \times \Gamma}$ and $u_0$.
}

We are concerned with the uniqueness and the stability of the determination of the coefficient $p.$
The results depend essentially on $t_0=0$, $0<t_0<T$ and
$t_0=T$.
\begin{itemize}
\item
Case $t_0=0$: No results in general.  We can refer to 
Klibanov \cite{K} in a case where a parabolic equation is considered
in $(0,T) \times \R^n$.  In \cite{K}, the uniqueness is proved
in an arbitrarily chosen ball $B \subset \R^n$ with boundary data
on $(0,T) \times \ppp B$. 
\item
Case $0<t_0<T$: Global Lipschitz stability as well as the global uniqueness.
See Imanuvilov and Yamamoto \cite{IY1}, and also Bellassoued and Yamamoto 
\cite{BY}.
\item
Case $t_0=T$: to our best knowledge no results are published.
\end{itemize}
In the case $t_0=0$, our problem is an inverse problem 
corresponding to a classical initial boundary value problem.
We note that in the case $t_0=T$, our observation ends at 
$T$, so that we have no extra boundary data
$u\vert_{\Gamma}$ for $t>T$. 
 
Our current article aims at two unsolved cases :
$t_0=0$ and $t_0=T$.
Our result in the case $t_0=0$ asserts the uniqueness in a more general
parabolic operator and subboundary and the proof is based on a similar
idea to \cite{K}.  The stability for $t_0=0$ is still totally an open problem.
On the other hand, in the case $t_0=T$, to the best knowledge of the authors, 
there are no existing results.  Our result for $t_0=T$ establishes not only 
the uniqueness 
but also the Lipschitz stability around any fixed coefficient $p$.

We emphasize that we consider the case of the boundary condition 
(1.2) of the Neumann type, not the Dirichlet boundary condition.  
For the global stability and the 
uniqueness by extra data on subboundary of $\ppp\OOO$, the case 
(1.2) requires more subtle arguments.  The arguments in the case of 
the Dirichelt boundary condition with extra Neumann data are easier, but 
owing to the required positivity assumption of spatial data $u_0$ in 
(1.3), we can not prove neither the global uniqueness nor the global 
stability without additional assumption that unknown coefficients 
are given in a neighborhood of $\ppp\OOO$.

For the case of $t_0=0$, our proof is via a corresponding inverse 
hyperbolic problem, and then a Carleman estimate for the Neumann 
boundary condition with additional Dirichlet data is indispensable 
and more works are required for establishing such a Carleman estimate 
(e.g., Imanuvilov \cite{Im1}).

We add that the boundary condition (1.2) means the thermal insulation, and
is quite important physically.    
\\

Henceforth we use the following notations: 
$Q := (0,T) \times \OOO$.

We arbitrarily fix  constants $M>0$, $\delta_0 > 0$, $\gamma_1, \gamma_2 \in  (0,1)$ and set
$$
\left\{ \begin{array}{rl}
& \mathcal{P} := \{ p\in C^{\gamma_1}(\ooo{\OOO});\,
\Vert p\Vert_{C^{\gamma_1}(\ooo{\OOO})} \le M\}, \\
& \mathcal{A} := \{ a\in C^{2+\gamma_2}(\ooo{\OOO}) ;\, 
\ppp_{\nu_A}a = 0 \,\, \mbox{on $\ppp\OOO$,}\\
& \quad \Vert a\Vert_{C^{2+\gamma_2}(\ooo{\OOO})} \le M, \quad 
a \ge \delta_0 > 0 \quad \mbox{on $\ooo{\OOO}$} \}.
\end{array}\right.
                                     \eqno{(1.7)}
$$

By $u_{p,a} = u_{p,a}(t,x)$ we denote the classical solution to 
$$
\left\{ \begin{array}{rl}
& \ppp_tu = Au - p(x)u \quad \mbox{in $Q$}, \\
& \ppp_{\nu_A}u = 0 \quad \mbox{on $(0,T) \times \ppp\OOO$},\\
& u(0,\cdot) = a \quad \mbox{in $\OOO$}.
\end{array}\right.
                                               \eqno{(1.8)}
$$
More precisely, we choose $\gamma > 0$ such that 
$$
0 < \gamma < \min\{ \gamma_1,\, \gamma_2\}<1.  
                                    \eqno{(1.9)}
$$
Then 
$$
u_{p,a} \in C^{1+\frac{\gamma}{2}, 2+\gamma}(\ooo{Q}) \quad
\mbox{and} \quad t\ppp_tu \in H^{1,2}(Q) \quad 
\mbox{for all $p\in \mathcal{P}$ and $a \in \mathcal{A}$}.
                                            \eqno{(1.10)}
$$
\\
{\bf Proof of (1.10).}
First, by Ladyzenskaja, Solonnikov and Ural'ceva \cite{LSU}, 
for $p\in \mathcal{P}$ and $a \in \mathcal{A}$, by noting (1.9) there exists 
a unique solution $u_{p,a} \in C^{1+\frac{\gamma}{2}, 2+\gamma}(\ooo{Q})$.
Moreover we see that there exists
a unique solution $v \in H^{1,2}(Q)$ to
$$
\left\{ \begin{array}{rl}
& \ppp_tv = Av - pv + \ppp_tu \quad \mbox{in $Q$},\\
& \ppp_{\nu_A}v = 0 \quad \mbox{on $(0,T)\times \ppp\OOO$}, \\
& v(0,\cdot) =0.
\end{array}\right.
$$
We can verify the unique existence of such $v$, similarly to 
Theorem 5 (pp.360-361) in Evans \cite{E} for example, where the homogeneous 
Dirichlet boundary condition is considered.
Then the uniqueness of the solution in $L^2(Q)$ yields that $v = t\ppp_tu$ in 
$Q$, which completes the proof of (1.10).
$\blacksquare$ 
\\

Now we state our main results.
\\
{\bf Theorem 1}
{\it
Let $\Gamma \subset \ppp\OOO$ be an arbitrarily chosen subboundary and 
$0 < \gamma < \min\{ \gamma_1, \, \gamma_2\}$.
In (1.3) we assume $t_0=T$.
There exists a constant $C_\gamma>0$, depending on $\mathcal{A}$ and 
$\mathcal{P}$ such that 
$$
\Vert p-q\Vert_{H^{\gamma}(\OOO)}
\le C_\gamma( (\Vert u_{p,a}- u_{q,b})(\cdot,T)\Vert_{H^{2+\gamma}(\OOO)}
+ \Vert u_{p,a} - u_{q,b}\Vert_{H^1((0,T)\times \Gamma)})
$$
for each $(p,a), (q,b) \in \mathcal{P} \times \mathcal{A}$.
}
\\

The global Lipschitz stability in the case of $t_0=T$ 
was an open problem.  However, Theorem 1 firstly provides the positive answer.
If spatial data are given at $t_0 \in (0,T)$, then we already have proved 
the global Lipschitz stability for the inverse parabolic problem 
(Imanuvilov and Yamamoto \cite{IY1}).
\\

Thanks to the holomorphicity of $u_{p,a}, u_{q,b}$ in $t>0$, 
fixing $\www{T} > T$ arbitrarily, by the three circle theorem 
(e.g., Chapter 10 in Cannon \cite{Ca}) we can obtain 
$$
\Vert u_{p,a} - u_{q,b}\Vert_{H^1((0,\www{T})\times \Gamma)}
\le C\Vert u_{p,a} - u_{q,b}\Vert_{H^1((0,T)\times \Gamma)}^{\mu}
$$
with $\mu \in (0,1)$ and $C>0$ within some boundedness assumptions.
Therefore in the case of data at the final time $T$, the essential 
contribution is the Lipschitz stability.
\\
 
Next we discuss the case of $t_0=0$, that is,
spatial data are given at the initial time moment in (1.1) - (1.3).
We consider a simpler parabolic equation.
$$
\left\{ \begin{array}{rl}
& \ppp_tu(t,x) = \Delta u (t,x) - p(x)u(t,x), \quad 0<t<T, \,
x\in\OOO, \\
& \ppp_{\nu_A}u = 0 \quad \mbox{on $(0,T) \times \ppp\OOO$}, \\
& u(0,\cdot) = u_0 \quad \mbox{in $\OOO$}.
\end{array}\right.
                                   \eqno{(1.11)}
$$
In this case, we note that $\ppp_{\nu_A}u = \nabla u\cdot \nu(x)$.
For simplicity, we assume that the spatial dimensions $n$ is smaller than or
equal to $3$: 
$$
n\in\{1,2,3\}.
$$
Furthermore we assume that $\OOO$ is locally convex in 
$\ppp\OOO \setminus {\ooo\Gamma_0}$, that is, for all 
$x\in \ppp\OOO \setminus 
{\ooo\Gamma_0}$, there exists a small ball $U_x$ such that 
$U_x \cap \OOO$ is convex.
We assume that $u_0\in H^3(\OOO)$ and $\ppp_{\nu_A}u_0 = 0$ on $\ppp\OOO$, and 
$p, q \in W^{1,\infty}(\OOO)$ for unknown coefficients.  
Replacing $u$ by $e^{C_0t}u$ with suitable constant $C_0$, without loss of 
generality, we can assume that
$$
p, q > 0 \quad \mbox{on $\ooo{\OOO}$}.           \eqno{(1.12)}
$$
By $u(p)$, we denote the solution to (1.11).

Arbitrarily choosing $x_0 \not\in \ooo{\OOO}$, we assume that 
a subboundary $\Gamma_0 \subset \ppp\OOO$ is a relatively open subset 
in $\R^{n-1}$ and
$$
\Gamma_0 \supset \{ x\in \ppp\OOO;\, (x-x_0)\cdot \nu(x) \ge 0\}.
                                      \eqno{(1.13)}
$$
Then we can state our second main result for the case of $t_0=0$.
\\
{\bf Theorem 2.}
\\
{\it 
We assume that $p,q \in W^{1,\infty}(\OOO)$ and 
$u_0 \in H^3(\OOO)$, $\ppp_{\nu_A}u_0 = 0$ on $\ppp\OOO$, and 
there exists a constant $\delta_0 > 0$ such that 
$$
\vert u_0(x) \vert \ge \delta_0 \quad \mbox{for all $x\in \ooo{\OOO}$}.
$$
Let $T>0$ be arbitrary.
Then $u(p) = u(q)$ on $(0,T) \times \Gamma_0$ implies
$p(x) = q(x)$ for $x\in \OOO$.
}

In Theorem 2, we can obtain the uniqueness for general parabolic 
operator as in (1.1), and we return to this issue in Section 6.

We prove Theorem 2 by reducing the inverse parabolic problem 
to an inverse problem for hyperbolic
equation and the argument relies on Klibanov \cite{K}.
The work \cite{K} assumes that $u(p), u(q)$ should be 
sufficiently smooth, and considers a parabolic equation 
in the whole space $(0,T)\times \R^n$ where unknown coefficients 
$p$, $q$ are given outside of a ball.

We note that since the proof is based on reduction of original problem to an inverse problem for a hyperbolic  equation, 
we need a Carleman estimate for hyperbolic equation  for the case of the Neumann boundary 
condition, and in Theorem 2 we should choose some large portion 
$\Gamma_0$ of the boundary such that (1.13) holds true, but $T>0$ can be arbitrarily small.
\\

Finally we discuss one type of inverse source problem.
Let 
$$
\ppp_ty = Ay - p(x)y + \mu(t)f(x), \quad 0<t<T, \, x\in \OOO, 
                                      \eqno{(1.14)}
$$
and
$$
y(t_0,\cdot) = 0 \quad \mbox{in $\OOO$}.
                                          \eqno{(1.15)}
$$
We assume that $\mu \not\equiv 0$ in $(0,T)$ is given.  Then 
\\  
{\bf Inverse source problem.}
\\
{\it 
Let $T>0$ be arbitrary and $0\le t_0 \le T$ be fixed.
Let $\Gamma \subset \ppp\OOO$ be an arbitrarily chosen subboundary.
Then determine $f$ in $\OOO$ by $y\vert_{(0,T) \times \Gamma}$ and
$\nabla y\vert_{(0,T)\times \Gamma}$.
}

We emphasize that we do not assume any boundary condition on the whole
lateral boundary $(0,T) \times \ppp\OOO$.  We are concerned with the 
uniqueness.
In the case of $0<t_0<T$, in view of the method based on  Carleman estimates, 
we can prove the uniqueness:
Assuming that $\mu\in C^1[0,T]$ and $\mu(t_0) \ne 0$, we can conclude 
that $y=\vert \nabla y\vert = 0$ on $(0,T) \times \Gamma$ and 
$y(t_0,\cdot) = 0$ in $\OOO$ imply $f=0$ in $\OOO$.
We omit the proof but we can refer for example to Yamamoto \cite{Ya2022}.

However, there are no published works in the case $t_0=0$ or $t_0=T$.
We establish the uniqueness in the case of $t_0=0$, which means the 
inverse source problem exactly corresponding to the initial value problem.
\\
{\bf Theorem 3.}
\\
{\it
Let $y, \ppp_ty \in H^{1,2}(Q)$ satisfy (1.14) and 
$$
y(0,\cdot) = 0 \quad \mbox{in $\OOO$}, 
$$
and $\mu \in C^1[0,T]$ and $\mu(0) \ne 0$. 
Then $y = \vert \nabla y\vert = 0$ on $(0,T)\times \Gamma$ 
implies $f=0$ in $\OOO$.}

We compare the uniqueness for the inverse source problem with the 
unique continuation which yields the uniqueness of solution without 
boundary condition on $\ppp\OOO$.
\\
{\bf Unique continuation.}
$$
\left\{ \begin{array}{rl}
& \ppp_tz = Az - p(x)z \quad \mbox{in $(0,T) \times \OOO$}, \\
& z = \vert \nabla z\vert = 0 \quad \mbox{on $(0,T)\times\Gamma$},
\end{array}\right.
$$
yield $z=0$ in $(0,T) \times \OOO$.
In particular, $z(0,\cdot) = 0$ in $\OOO$.
\\
{\bf Uniqueness in the inverse source problem.}
$$
\left\{ \begin{array}{rl}
& \ppp_ty = Ay - p(x)y + \mu(t)f(x)
\quad \mbox{in $(0,T) \times \OOO$}, \\
& y = \vert \nabla y\vert = 0 \quad \mbox{on $(0,T)\times\Gamma$},\\
& y(0,\cdot)=0 \quad \mbox{in $\OOO$},
\end{array}\right.
$$
yield $f=0$ in $\OOO$, and so $y=0$ in $(0,T) \times \OOO$.
\\

Comparing with the unique continuation, in our inverse source problem, one 
more spatial function $f$ is unknown, and for compensating for the uniqueness,
we need more spatial data at $t=0$, but never require additional boundary data.
Thus Theorem 3 provides the uniqueness result with the minimum data.

As for conditional stability results for the case of $t_0=0$ with 
boundary condition on 
$(0,T) \times \ppp\OOO$, see Yamamoto \cite{Ya1992} for example.
\\

This article is composed of five sections.
In Section 2, we prove Theorem 1 by a priori estimates of 
$u(p) - u(q)$ involving not initial values and the compact-uniqueness 
argument.
Sections 3 and 4 are devoted to the proofs of
Theorem 2 and 3 respectively.
Section 5 gives concluding remarks.

\section{Proof of Theorem 1}

We recall that $\Gamma \subset \ppp\OOO$ is an arbitrarily chosen 
non-empty subboundary.  We set 
$$
\left\{ \begin{array}{rl}
& Q:= (0,T) \times \OOO, \quad 
Q_1:= (\frac{T}{2}, \,T) \times \OOO, \quad
Q_2:= (\frac{T}{4}, \,T) \times \OOO, \\
& Q_3:= (\frac{T}{8}, \,T) \times \OOO, \quad
  \Sigma := (0,T) \times \Gamma.
\end{array}\right.
                                   \eqno{(2.1)}
$$
\\
{\bf Proof.}
We divide the proof into five steps.
The key of the proof is:
\\
(i) estimate of $\Vert p-q\Vert_{H^{\gamma}(\OOO)}$ by data of 
$u_{p,a}-u_{q,b}$ and $\Vert p-q\Vert_{L^2(\OOO)}$, which is given by 
(2.18) below.
\\
(ii) application of the compact-uniqueness argument in order to 
eliminates $\Vert p-q\Vert_{L^2(\OOO)}$ in (2.18).
\\

{\bf First Step.}
In order to formulate a  Carleman estimate we introduce a weight function for 
such a Carleman estimate by the following lemma:
\\
{\bf Lemma 2.1 (\cite{FI})}
\\
{\it 
Let $\Gamma_0$ be an arbitrarily fixed subboundary of $\partial\Omega$
such that $\ooo{\Gamma_0} \subset \Gamma$.
Then there exists a function $\psi \in C^2(\overline{\Omega})$ such that
$$
\psi(x)>0 \quad \mbox{for all $x\in\Omega$}, \quad
\psi|_{\partial\Omega\setminus\Gamma_0} = 0, \quad |\nabla \psi(x)| > 0 \quad
\mbox{for all $x\in \ooo{\Omega}$}.
$$
}

Using the function $\psi$, we define 
%\begin{equation}\label{DD1.2} 
$$
\alpha(t, x) = \frac{e^{\lambda\psi(x)}
- e^{2\lambda ||\psi||_{C(\overline{\Omega})}}}{t(T-t)}, \quad
\varphi(t, x) = \frac{e^{\lambda \psi(x)}}{t(T-t)}.
$$
Consider the boundary value problem:
%\begin{equation}\label{D1.4}
$$
Pz(t,x):= \ppp_tz - Az + p(x)z =g\quad \mbox{in $Q$}, 
\quad \partial_{\nu_A}z = 0 \quad \mbox{on $(0,T) \times \ppp\OOO$}.
$$
%\end{equation}
We are ready to state a Carleman estimate.
\\
%\begin{lemma}\label{D1.2}
{\bf Lemma 2.2 (\cite{FI})}
\\
{\it Let (1.4) and (1.5) be fulfilled.
Then there exists a constant $\widehat{\lambda} > 0$ such that for 
an arbitrary $\lambda \ge \widehat{\lambda}$, there exist constants
${\tau}_0(\lambda)>0$ and $C_1>0$ such that 
\begin{align*}
& \int_Q \left(\frac{1}{{\tau}\varphi} (\vert\partial_t z\vert^2 
+ \sum_{j,k=1}^n\vert \partial_j\ppp_kz\vert^2) 
+ {\tau}\varphi |\nabla z|^2 + {\tau}^3 \varphi^3 z^2\right)
e^{2{\tau}\alpha}   dx\,dt\\
\le& C_1\left( \int_Q |Pz|^2 e^{2{\tau}\alpha}  dxdt +
\int_{(0,T)\times \Gamma} \{\tau \varphi( \vert \nabla z\vert^2 
+ \vert\partial_t z\vert^2) 
+ (\tau\varphi)^3 \vert z\vert^2\} e^{2{\tau}\alpha}  d\Sigma\right)
\end{align*}
for each $\tau \ge \tau_0(\lambda)$ and all $z,t\partial_t z \in H^{1,2}(Q)$ satisfying 
$\ppp_{\nu_A}z = 0$ on $(0,T) \times \ppp\OOO$.
}
\\

The constant $C_1 > 0$ can be taken 
independently of choices of $p\in \PPPP$.  Indeed $C_1$ depends on
$M_0$, provided that $\Vert p\Vert_{C(\ooo{\OOO})} \le M_0$.
\\
We recall that $u_{p,a}$ is the solution to (1.8) and  
$0<\gamma < 1$ satisfies (1.9).
\\

Furthermore, in addition to (1.10), we can prove 
\\
{\bf Lemma 2.3.}
\\
{\it  
For any $\gamma_*$ satisfying 
$0< \gamma_* < \min\{ \gamma_1, \, \gamma_2\}$, there 
exists a constant $C>0$ such that 
$$
\Vert u_{p,a}\Vert_{C^{1+\frac{\gamma_*}{2},2+\gamma_*}(\ooo{Q})} 
\le C, \quad 
\Vert t \ppp_tu_{p,a}\Vert_{H^{1,2}(Q)} \le C \quad 
\mbox{for all $p\in \PPPP$ and $a\in \AAAA$}.
$$
}
Indeed, the estimates of 
$\Vert u_{p,a}\Vert_{C^{1+\frac{\gamma^*}{2},2+\gamma^*}(\ooo{Q})}$ and 
$\Vert t\ppp_tu_{p,a}\Vert_{H^{1,2}(Q)}$ follow respectively, for example
from Theorem 5.3 (pp.320-321) in \cite{LSU} and 
Theorem 5 (pp. 360-361) in \cite{E}.
\\

Moreover,
\\
{\bf Lemma 2.4.}
{\it There exist a strictly positive constant $c_0(\delta_0)$ such that 
$$
u_{q,b}(T, \cdot) \ge c_0(\delta_0) \quad \mbox{on $\ooo{\OOO}$}.          
$$ 
for all $b\in \AAAA$ and $q \in \PPPP$}

The proof is done by the positivity of the fundamental solution 
(e.g., It\^o \cite{It}) in terms of $b \ge \delta_0$ on $\ooo{\OOO}$
for $b \in \AAAA$ and $q\in \PPPP$.

{\bf Second Step.}
For any $\mbox{\bf v} = (v_1,..., v_n), \mbox{\bf w}
= (w_1,..., w_n)\in \Bbb R^n$, we set 
%\begin{equation}\label{pioner}
$$
a(x,\mbox{\bf v},\mbox{\bf w})=\sum_{j,k=1}^n a_{jk}(x)v_jw_k.
$$
Let $u \in H^{1,2}(Q)$ satisfy $t\ppp_tu \in H^{1,2}(Q)$ and 
$$
\left\{ \begin{array}{rl}
& \ppp_tu = Au - p(x)u + S(t,x) \quad \mbox{in $Q$}, \\
& \ppp_{\nu_A}u = 0 \quad \mbox{on $(0,T)\times \ppp\OOO$}
\end{array}\right.
                                               \eqno{(2.2)}
$$
for $S \in L^2(Q)$ and $p\in \mathcal{P}$.  In particular, we note that 
$\Vert p\Vert_{C(\ooo{\OOO})} \le M$, where the constant $M>0$ is given in 
(1.7).

Thanks to Lemma 2.2, we can readily prove
%\begin{equation}\label{sprout}
$$
\Vert u\Vert_{L^2(\frac{T}{4},\frac{3T}{4};H^1(\Omega))}
\le C(\Vert S\Vert_{L^2(Q_3)} + \Vert u\Vert_{H^1(\Sigma)}). 
                                                         \eqno{(2.3)}
$$
%\end{equation}
Here we recall (2.1).
In this step, we will further prove
%\begin{equation}\label{Ssprout}
$$
\Vert u\Vert_{L^2(\frac{T}{4},T ;H^1(\Omega))}\le C(\Vert S\Vert_{L^2(Q_3)}
+ \Vert u\Vert_{H^1(\Sigma)}).              \eqno{(2.4)}
%\end{equation}
$$
{\bf Proof of (2.4).}
We choose $\mu\in C^\infty[\frac{T}{4},\frac{3T}{4}]$ satisfying 
$\mu\left(\frac{T}{4}\right)=0$ and $\mu\left(\frac{3T}{4}\right)=1$.
We take the scalar products of equation (2.2) with $\mu(t) u$ 
in $L^2((\frac T4,\frac{3T}{4})\times \Omega)$.
Integrating by parts, we have
\begin{align*}
& \int_{\frac T4}^{\frac {3T}{4}}\int_\Omega\biggl( 
\frac 12 \partial_t(\mu(t)u^2)-\frac{\partial_t\mu(t)}{2}u^2
+ \mu(t) a(x,\nabla u,\nabla u)  
+  \left( \sum_{j=1}^n b_j(x)\partial_ju+p(x) u\right)\mu u\biggr)dxdt\\
= &\int_{\frac T4}^{\frac {3T}{4}}\int_\Omega\mu(t)uS\, dxdt.
\end{align*}
This equality implies
\begin{align*}
& \left\Vert u\left(\frac {3T}{4},\cdot\right)\right\Vert^2_{L^2(\Omega)}
= 2\int_{\frac T4}^{\frac {3T}4}\int_\Omega\biggl( 
\frac{\partial_t\mu(t)}{2}u^2-2\mu(t) a(x,\nabla u,\nabla u)\\
-2& \left(\sum_{j=1}^n b_j(x)\partial_ju+p(x) u\right)\mu u\biggr)dxdt
+ 2\int_{\frac T4}^{\frac {3T}4}\int_\Omega\mu(t)uS\,dxdt
\end{align*}
$$
\le C(\Vert u\Vert^2_{L^2(\frac{T}{4},\frac{3T}{4};H^1(\Omega))}
+ \Vert S\Vert^2_{L^2(Q_3)}).                      \eqno{(2.5)}
$$
Applying (2.3) to the first term on the right-hand side of (2.5),
we have 
%\begin{equation}\label{samson}
$$
\left\Vert u\left(\frac {3T}{4},\cdot\right)\right\Vert_{L^2(\Omega)}
\le  C(\Vert S\Vert_{L^2(Q_3)}+\Vert u\Vert_{H^1(\Sigma)}).
                                                 \eqno{(2.6)}
$$
%\end{equation}
In view of (2.6), 
the standard a priori estimate for the initial value problem with 
initial value at $\frac{3}{4}T$ (e.g., Section 1 of Chapter 7 in 
Evans \cite{E}), yields 
$$
\Vert u\Vert_{L^2(\frac{3T}{4},T;H^1(\Omega))}
\le C\left( \left\Vert u(\frac {3T}{4},\cdot)\right\Vert_{L^2(\Omega)}
+ \Vert S\Vert_{L^2(Q_3)}\right)
\le C(\Vert S\Vert_{L^2(Q_3)}+\Vert u\Vert_{H^1(\Sigma)}).
$$
Combining this estimate with (2.3), we complete the proof of (2.4).
$\blacksquare$
\\
{\bf Third Step.}
We assume that $u \in H^{1,2}(Q)$ satisfies $t\ppp_tu \in H^{1,2}(Q)$ and (2.2).
In this step, we first prove
%\begin{equation}\label{sprout1}
$$
\Vert u\Vert_{L^2(\frac T2 ,T;H^2(\Omega))}
+ \Vert \partial_tu\Vert_{L^2(\frac T2 ,T;L^2(\Omega))}
\le C(\Vert S\Vert_{L^2(Q_3)} + \Vert u\Vert_{H^1(\Sigma)}).
                                                      \eqno{(2.7)}
$$
%\end{equation}
{\bf Proof of (2.7).}
Since $a_{jk} = a_{kj}$ are independent of $t$, in view of 
integration by parts and $\ppp_{\nu_A}u = 0$ on $\ppp\OOO$, we see 
\begin{align*}
& -\int_{\OOO} \sum_{j,k=1}^n \ppp_j(a_{jk}\ppp_ku)\ppp_tu dx 
= \sumjk a_{jk}(\ppp_ku)\ppp_j\ppp_t u dx\\
=& \int_{\OOO} \sum_{k>j} a_{jk}((\ppp_ku)\ppp_j\ppp_tu
+ (\ppp_ju)\ppp_k\ppp_tu) dx 
+ \int_{\OOO} \sum_{j=1}^n a_{jj}(\ppp_ju)\ppp_j\ppp_tu dx \\
=& \int_{\OOO} \sum_{k>j} a_{jk} \ppp_t((\ppp_ju)\ppp_ku)dx 
+ \frac{1}{2} \int_{\OOO} \sum_{j=1}^n a_{jj}\ppp_t((\ppp_ju)^2) dx \\
=& \frac{1}{2} \int_{\OOO} \ppp_t\left( 
\sumjk a_{jk}(\ppp_ju)\ppp_ku \right) dx.
\end{align*}
Taking the scalar products of equation (2.2) by 
$\left(t - \frac{T}{4}\right)\partial_t u$ in $L^2((\frac T4,T)\times \Omega)$
and integrating by parts, we obtain
\begin{align*}
& \int_{\frac T4}^T\int_\Omega\biggl((t-\frac T4)(\partial_tu)^2
+\frac 12\partial_t\left( (t-\frac{T}{4})
a(x,\nabla u,\nabla u)\right)
-\frac 12 a(x,\nabla u,\nabla u)\\
+ &\left(\sum_{j=1}^n b_j(x)\partial_ju+p(x) u\right)
(t-\frac T4)\partial_tu\biggr)dxdt\\
=& \int_{\frac{T}{4}}^T\int_\Omega ((t-\frac{T}{4})(\partial_tu)^2
-\frac{1}{2} a(x,\nabla u,\nabla u)
+ \left( \sum_{j=1}^n b_j(x)\partial_ju+p(x) u\right)(t-\frac T4)
\partial_tu ) dxdt\\
+& \frac{3T}{8} \int_\Omega a(x,\nabla u(T,x),\nabla u(T,x))dx\\
= & \int_{\frac{T}{4}}^T \int_\Omega(t-\frac T4)(\partial_tu)S\, dxdt.
\end{align*}
Hence, by $a(x,\nabla u, \nabla u) \ge 0$, we have
\begin{align*}
& \int_{\frac T4}^T\int_\Omega(t-\frac T4)(\partial_tu)^2dxdt
\le \int_{\frac T4}^T \int_{\OOO} \left( \frac 12 a(x,\nabla u,\nabla u)
+ \left(\sum_{j=1}^n b_j(x)\partial_ju+p(x) u\right)(t-\frac T4)\partial_tu
\right)dxdt\\
+ & \int_{\frac T4}^T\int_\Omega(t-\frac T4)(\partial_tu)S\, dxdt.
\end{align*}

Here for any $\ep \in (0,1)$, we can choose a constant $C_{\ep} > 0$ such that 
$$
\left\vert (t-\frac{T}{4}) (\ppp_tu)S \right\vert
\le \ep (t-\frac{T}{4})\vert \ppp_tu\vert^2
+ C_{\ep}(t-\frac{T}{4})\vert S\vert^2
\le \ep(t-\frac{T}{4})\vert \ppp_tu\vert^2
+ C_{\ep}\frac{3T}{4}\vert S\vert^2 
$$
in $\left( \frac{T}{4}, T\right) \times \OOO$.
Consequently, 
$$
\int_{\frac T4}^T\int_\Omega(t-\frac T4)(\partial_tu)^2dxdt
\le C_{\ep}\left(\Vert u\Vert^2_{L^2(\frac{T}{4},T;H^1(\Omega))}
+ \int_{\frac T4}^T\int_\Omega\vert S\vert^2dxdt\right)
+ \ep\int_{\frac T4}^T\int_\Omega(t-\frac T4)(\partial_tu)^2dxdt.
$$
Fixing $\ep \in (0,1)$ and applying (2.4), we obtain
%\begin{equation}\label{villa}
$$
\Vert \partial_tu\Vert^2_{L^2(\frac T2 ,T;L^2(\Omega))}
\le \int_{\frac T4}^T\int_\Omega(t-\frac T4)(\partial_tu)^2dxdt
$$
$$
\le C(\Vert u\Vert^2_{L^2(\frac{T}{4},T;H^1(\Omega))}
+ \Vert S\Vert^2_{L^2(Q_3)})
\le C(\Vert S\Vert^2_{L^2(Q_3)}+\Vert u\Vert^2_{H^1(\Sigma)}).
                                 \eqno{(2.8)}
$$
%\end{equation}
In order to estimate the second space derivatives of $u$,
we will write the equation as
$$
-\sum_{j,k=1}^n \partial_{x_j}(a_{jk}(x)\partial_ku)
= -\partial_tu-\sum_{j=1}^n b_j(x)\partial_ju-p(x) u + S(t,x),
\quad x\in \OOO, \, 0<t<T
$$
with $\ppp_{\nu_A}u = 0$ on $(0,T) \times \ppp\OOO$.
In particular, we have
\begin{align*}
& \left\Vert -\sum_{j,k=1}^n \partial_j(a_{jk}(x)\partial_ku(t,\cdot))
\right\Vert_{L^2(\OOO)}\\
\le& C(\Vert \ppp_tu(t,\cdot)\Vert_{L^2(\OOO)}
+ \Vert \nabla u(t,\cdot)\Vert_{L^2(\OOO)}
+ M\Vert u(t,\cdot)\Vert_{L^2(\OOO)}
+ \Vert S(t,\cdot)\Vert_{L^2(\OOO)}), \quad \frac{T}{4} \le t \le T.
\end{align*}

Let us fix $t\in (\frac T2,T)$.
From the a priori estimate for the Neumann problem for the second order 
elliptic operator, we obtain
$$
\Vert u(t,\cdot)\Vert^2_{H^2(\Omega)}
\le C(\Vert\partial_tu(t,\cdot)\Vert^2_{L^2(\Omega)}
+ \Vert u(t,\cdot)\Vert^2_{H^1(\Omega)}+\Vert S(t,\cdot)\Vert^2_{L^2(\Omega)}).
$$

Integrating this inequality over the time interval $(\frac T2 ,T)$,
we obtain
%\begin{equation}\label{villa1}
$$
\Vert u(t,\cdot)\Vert^2_{L^2(\frac T2,T;H^2(\Omega))}
\le C(\Vert\partial_tu(t,\cdot)\Vert^2_{L^2(\frac T2,T;L^2(\Omega))}
+ \Vert u(t,\cdot)\Vert^2_{L^2(\frac T2,T;H^1(\Omega))}
+ \Vert S\Vert^2_{L^2(Q_3)}).
                                                \eqno{(2.9)}
$$
%\end{equation}
Using estimates (2.8) and (2.4) to estimate the right hand side of (2.9),
we reach (2.7).
$\blacksquare$
\\

{\bf Fourth Step.}
Let $u \in H^{1,2}(Q)$ satisfy $t\ppp_tu \in H^{1,2}(Q)$ and (2.2).

By $t\ppp_tu \in H^{1,2}(Q)$, we differentiate equation (2.2) with 
respect to $t$:
%\begin{equation}\label{VSsok1}
$$
\ppp_t(\ppp_tu) = A(\ppp_tu) - p(x)\ppp_tu + \ppp_tS(t,x)f(x)
\quad \mbox{in $Q$},
                                         \eqno{(2.10)}
$$
$$
\partial_{\nu_A} \ppp_t u = 0 \quad \mbox{on $(0,T) \times \ppp\OOO$}
                                           \eqno{(2.11)}
$$
and
$$
(\ppp_tu)(T,\cdot) = (A-p(x))u(T,\cdot) + S(T,\cdot).
                                           \eqno{(2.12)}
$$
We recall that $\gamma \in (0,1)$ satisfies (1.9).
Now we prove
$$
\Vert S(T, \cdot)\Vert_{H^{\gamma}(\Omega)}
\le C(\Vert \ppp_tS\Vert_{L^2(Q_1)} + \Vert S\Vert_{L^2(Q_3)}
+ \Vert u\Vert_{H^1(\Sigma)}
+ \Vert u(T,\cdot)\Vert_{H^{2+\gamma}(\OOO)}). 
                                            \eqno{(2.13)}
$$
\\
{\bf Proof of (2.13).}
We set $\widetilde u = (t-\frac T2)\partial_tu$.  Then
%\begin{equation}\label{VVSsok1}
$$
\left\{ \begin{array}{rl}
& \ppp_t\www{u} - A\www{u} + p(x)\www{u} 
= (t-\frac T2)\partial_tS(t,x) + \partial_t u \quad 
\mbox{in $\left(\frac {T}{2},T\right)\times \Omega$},\\
& \partial_{\nu_A}\widetilde u\vert_{(\frac {T}{2},T)\times \partial\Omega}=0, 
\quad \widetilde u(\frac T 2,\cdot)=0.
\end{array}\right.
                                               \eqno{(2.14)}
$$
%\end{equation}
Combining the classical a priori estimate yields
$$
\www{u} \in L^2( \frac{T}{2},T; H^2(\OOO))
\cap H^1( \frac{T}{2},T; L^2(\OOO))
$$
(e.g., Theorem 5 (p.360) in Evans \cite{E}) with an interpolation 
property (Theorem 3.1 (p.19) in Lions and Magenes \cite{LM}), we can see
the continuous embedding 
$$
L^2( \frac{T}{2},T; H^2(\OOO))
\cap H^1( \frac{T}{2},T; L^2(\OOO))
\subset C( [\frac{T}{2},T]; H^1(\OOO)).
$$
We can directly verify this embedding without the interpolation 
property, but we do not use such verification.

Therefore,  
$$
\Vert \widetilde u\Vert_{C([\frac{T}{2},T]; H^1(\OOO))}
\le C\left\Vert (t-\frac T2)\partial_tS
+ \partial_t u\right\Vert_{L^2(Q_1)}
$$
$$
\le C(\Vert \ppp_tS\Vert_{L^2(Q_1)}
+ \Vert \partial_t u\Vert_{L^2(Q_1)}).
                                                   \eqno{(2.15)}
$$
Applying (2.7) in order to estimate the last term on the right-hand side 
of (2.15), we have
%\begin{equation}\label{sotka1}
$$
\left\Vert  (t-\frac T2)\partial_tu\right\Vert
_{C\left( \left[\frac{T}{2},T\right]; H^1(\OOO)\right)}
\le C(\Vert \ppp_tS\Vert_{L^2(Q_1)} + \Vert S\Vert_{L^2(Q_3)}
+ \Vert u\Vert_{H^1(\Sigma)}). 
                                                 \eqno{(2.16)}
$$
%\end{equation}

From (2.12), we have
$$
\Vert S(T,\cdot)\Vert_{H^\gamma(\OOO)}
\le C(\Vert \ppp_tu(T,\cdot)\Vert_{H^1(\OOO)}
+ \Vert (A-p)u(\cdot,T)\Vert_{H^\gamma(\OOO)})
$$
$$\le C(\Vert \ppp_tu(T,\cdot)\Vert_{H^1(\OOO)}
+ (1+\Vert p\Vert_{C^{\gamma_1}(\overline\Omega)})\Vert u(\cdot,T)\Vert_{H^{2+\gamma}(\OOO)}) .
$$
Thus we proved that (2.13) holds  
provided that $u\in H^{1,2}(Q)$ satisfies $t\ppp_tu \in H^{1,2}(Q)$ and (2.2).
$\blacksquare$
\\

{\bf Fifth Step.}
Let $z = z_{p,q,a,b}=u_{p,a}-u_{q,b} \in H^{1,2}(Q)$ and $f:=q-p
\in L^2(\OOO)$.  Then, by Lemma 2.3, we see that  
$t\ppp_tz \in H^{1,2}(Q)$ and
$$
\left\{ \begin{array}{rl} 
& \ppp_tz = Az - p(x)z + f(x)u_{q,b}(t,x) \quad \mbox{in $Q$},\\
& \partial_{\nu_A} z = 0 \quad \mbox{on $(0,T) \times \ppp\OOO$}, \\
& z(0,\cdot) = a-b \quad \mbox{in $\OOO$}.
\end{array}\right.
                                           \eqno{(2.17)}
$$
Set $S(t,x) := f(x)u_{q,b}(t,x)$ in $Q$.  By (2.1) and Lemmata 2.3 and 2.4 
we have
$$
\Vert f\Vert_{H^{\gamma}(\OOO)} \le \left\Vert
\frac{S(T,\cdot)}{u_{q,b}(T,\cdot)}\right\Vert_{H^{\gamma}(\OOO)}
\le C\Vert S(T,\cdot)\Vert_{H^{\gamma}(\OOO)}.
$$
Therefore, applying (2.13) to (2.17), we can find a constant 
$\widetilde C>0$ such that 
$$
\Vert p-q\Vert_{H^{\gamma}(\OOO)}
\le \widetilde C(\Vert p-q\Vert_{L^2(\OOO)} 
+ \Vert z_{p,q,a,b}(T,\cdot)\Vert_{H^{2+\gamma}(\OOO)}
+ \Vert z_{p,q,a,b}\Vert_{H^1(\Sigma)})                   \eqno{(2.18)}
$$
for all $p, q \in \PPPP$ and $a,b \in \AAAA$.
\\

Now we will prove that there exists a constant $C>0$ such that 
$$
\Vert p-q\Vert_{H^{\gamma}(\OOO)} 
\le C(\Vert z_{p,q,a,b}(T,\cdot)\Vert_{H^{2+\gamma}(\OOO)}
+ \Vert z_{p,q,a,b}\Vert_{H^1(\Sigma)})                    \eqno{(2.19)}
$$
for all $p,q \in \PPPP$ and $a,b \in \AAAA$.

We start the proof of (2.19).  Suppose that (2.19) does not hold. 
Then for any $m\in \Bbb N$, there exist functions $p_m,q_m \in \PPPP$ 
and $a_m,b_m \in \AAAA$ such that
$$
\Vert p_m-q_m\Vert_{H^{\gamma}(\OOO)} 
> m(\Vert z_{p_m,q_m,a_m,b_m}(T,\cdot)\Vert_{H^{2+\gamma}(\OOO)}
+ \Vert z_{p_m,q_m,a_m,b_m}\Vert_{H^1(\Sigma)}).
$$
Then, by the definition (1.8) of $u_{p_m,a_m}, u_{q_m,b_m}$, we can extend 
$u_{p_m,a_m}, u_{q_m,b_m}$ to $t\in (0, +\infty)$, so that 
$$
\left\{ \begin{array}{rl}
& \ppp_tu_{p_m,a_m} = Au_{p_m,a_m} - p_mu_{p_m,a_m}, \quad (t,x) \in 
(0, +\infty) \times \OOO, \\
& \partial_{\nu_A}u_{p_m,a_m}=0 \quad \mbox{on $(0,+\infty)\times 
\ppp\OOO$}, \\
& u_{p_m,a_m}(0,\cdot) = a_m \quad \mbox{in $\OOO$}
\end{array}\right.
                                                    \eqno{(2.20)}
$$
and
$$
\left\{ \begin{array}{rl}
& \ppp_tu_{q_m,b_m} = Au_{q_m,b_m} - q_mu_{q_m,b_m}, \quad (t,x) \in 
(0, +\infty) \times \OOO, \\
& \partial_{\nu_A}u_{q_m,b_m}=0 \quad \mbox{on $(0,+\infty)\times \ppp\OOO$},\\
& u_{q_m,b_m}(0,\cdot) = b_m \quad \mbox{in $\OOO$}.
\end{array}\right.
                                    \eqno{(2.21)}
$$
Moreover for any $\www{T}>0$, by Lemma 2.3, we see that 
$$
u_{p_m,a_m}, \, u_{q_m,b_m} \in C^{1+\frac{\gamma}{2}, 2+\gamma}
([0,\www{T}]\times \ooo{\OOO})
$$
and there exists a constant $C(\www{T}) > 0$ such that 
$$
\Vert u_{p_m,a_m}\Vert_{C^{1+\frac{\gamma}{2}, 2+\gamma}
([0,\www{T}]\times \ooo{\OOO})}, \,\,
\Vert u_{q_m,b_m}\Vert_{C^{1+\frac{\gamma}{2}, 2+\gamma}
([0,\www{T}]\times \ooo{\OOO})}\, \le C(\www{T}).
$$
Moreover, $z_{p_m,q_m,a_m,b_m}=u_{p_m,a_m}-u_{q_m,b_m}\in  H^{1,2}(Q)$ and
$$
\frac{\Vert z_{p_m,q_m,a_m,b_m}(T,\cdot)\Vert_{H^{2+\gamma}(\Omega)}
+ \Vert z_{p_m,q_m,a_m,b_m}\Vert_{H^1(\Sigma)}}
{\Vert p_m-q_m\Vert_{H^{\gamma}(\Omega)}}\, \rightarrow 0
\quad\mbox{as $m\rightarrow +\infty$}.
                                       \eqno{(2.22)}
$$
We set  
$$
w_m= \frac{z_{p_m,q_m,a_m,b_m}}{\Vert p_m-q_m\Vert_{H^{\gamma}(\Omega)}},
\quad r_m = \frac{q_m-p_m}{\Vert p_m-q_m\Vert_{H^{\gamma}(\Omega)}}.
$$
From (2.20) and (2.21) it follows that 
$w_m$ satisfies  the equations
$$
\left\{ \begin{array}{rl} 
& \ppp_tw_m = A w_m - p_m(x)w_m + r_m(x)u_{q_m,b_m}(t,x) \quad 
\mbox{in $(0,+\infty) \times \OOO$},\\
& \partial_{\nu_A} w_m = 0 \quad \mbox{on $(0,+\infty) \times \ppp\OOO$}, \\
\end{array}\right.
                                           \eqno{(2.23)}
$$
By Lemma 2.3 where we set $\gamma_* = \gamma$, 
the sequence $u_{q_m,b_m}$ is bounded in 
$C^{1+\frac{\gamma}{2},2+\gamma}(\ooo{Q})$. 
Therefore, since the sequence $p_m$ is bounded in $C^{\gamma_1}(\ooo\Omega)$
and $\Vert r_m\Vert_{L^2(\OOO)} \le 1$ for $m\in\N$,
applying Lemma 2.2 to (2.23), for any $\varepsilon>0$ we obtain 
$$
\Vert w_m\Vert_{H^{1,2}((\ep,T-\varepsilon)\times\Omega)}
\le C_\varepsilon.
$$
We fix $\www{T}>0$ arbitrarily.
The a priori estimate for parabolic equations implies
$$
\Vert w_m\Vert_{H^{1,2}((3\ep,\www{T})\times\Omega)}\le C_\varepsilon.
                                               \eqno{(2.24)}
$$
Indeed, we observe that $v_m=(t-2\varepsilon)\partial_t w_m$ solves the 
initial boundary value problem
$$
\left\{ \begin{array}{rl} 
& \ppp_tv_m = Av_m - p_m(x)v_m +\partial_t w_m 
+ r_m(t-2\ep)\partial_tu_{q_m,b_m}(t,x)
 \quad \mbox{in $(2\ep,\www{T})\times \Omega$},\\
& \partial_{\nu_A} v_m = 0 \quad \mbox{on}\, (0,\www{T}) \times \ppp\OOO,
\quad  v_m(2\varepsilon,\cdot) = 0.
\end{array}\right.                                      
                                           \eqno{(2.25)}
$$
It is known that the solution $v_m$ to (2.25) belongs to the space 
$H^{1,2}((2\varepsilon,\widetilde{T})\times \Omega)$ for all 
$\widetilde T>2\varepsilon$ (e.g., Theorem 5 (pp.360-361) in \cite{E}).
Thus (2.24) is verified.
$\blacksquare$
\\

By (2.24), we have 
$$
\Vert \ppp_jw_m\Vert_{L^2((3\ep,\www{T}) \times \OOO)}, \quad
\Vert \ppp_j\ppp_kw_m\Vert_{L^2((3\ep,\www{T}) \times \OOO)}, \quad
\Vert \ppp_tw_m\Vert_{L^2((3\ep,\www{T}) \times \OOO)} \le C_{\ep}
$$
for $1\le j,k\le n$.  Therefore there exist $w$, $W_j$, $W_{jk}$,
$W_0 \in L^2((3\ep,\www{T}) \times \OOO)$ such that we can extract 
subsequences, denoted by the same notations, to have
$w_m \longrightarrow w$, $\ppp_jw_m \longrightarrow W_j$,
$\ppp_j\ppp_kw_m \longrightarrow W_{jk}$, $\ppp_tw_m \longrightarrow W_0$
weakly in $L^2((3\ep,\www{T}) \times \OOO)$ as $m\to +\infty$ 
for $1\le j,k \le n$. Noting that 
$$
_{\mathcal{D}'}\langle \ppp_jw_m,\, \, \va\rangle_{\DDD}
= \, -\,_{\DDD'}\langle w_m,\, \, \ppp_j\va \rangle_{\DDD}\,
\longrightarrow \, -\,_{\DDD'}\langle w, \,\, \ppp_j\va\rangle
\quad \mbox{as $m\to \infty$}
$$
for all $\va \in \DDD:= C^{\infty}_0((3\ep,\www{T}) \times \OOO)$
and $\ppp_jw_m \longrightarrow W_j$ in $\DDD'$, etc., we can easily 
verify that $W_j=\ppp_jw$, $W_{jk} = \ppp_j\ppp_kw$,
$W_0 = \ppp_tw$ for $1\le j,k \le n$.  Here $\DDD'$ denotes the space of 
distributions in $(3\ep,\www{T}) \times \OOO$.
Hence,
$$
\left\{ \begin{array}{rl}
& w_m \longrightarrow w, \quad
  \ppp_jw_m \longrightarrow \ppp_jw, \quad
  \ppp_j\ppp_kw_m \longrightarrow \ppp_j\ppp_kw, 
\quad 1\le j,k \le n,                    \\
& \ppp_tw_m \longrightarrow \ppp_tw \quad \mbox{weakly in 
   $L^2((3\ep,\www{T}) \times \OOO)$ as $m\to +\infty$}.
\end{array}\right.
                                             \eqno{(2.26)}
$$
We arbitrarily fix $\gamma_*$ satisfying 
$0< \gamma_* < \gamma$.
Moreover, the Ascoli-Arzel\`a theorem yields that the embeddings
$C^{\gamma_1}(\ooo{\OOO}) \subset C^{\gamma_*}(\ooo{\OOO})$ and
$C^{1+\frac{\gamma}{2}, 2+\gamma}([0,\www{T}] \times \ooo{\OOO})
\subset C^{1+\frac{\gamma_*}{2}, 2+\gamma_*}([0,\www{T}] \times \ooo{\OOO})$
are compact by $\gamma_* < \gamma$.
Taking subsequences, denoted by the same notations, we obtain 
$$
p_m \longrightarrow p, \, q_m \longrightarrow q 
\quad \mbox{in $C^{\gamma_*}(\ooo{\OOO})$}, \quad
u_{q_m,b_m} \longrightarrow u_{q,b} \quad
\mbox{in $C^{1+\frac{\gamma_*}{2}, 2+\gamma_*}([0,\www{T}] \times 
\ooo{\OOO})$}
                                                         \eqno{(2.27)}
$$ 
as $m\to +\infty$ 
by $p_m, q_m\in \mathcal{P}$, $b_m \in \mathcal{A}$, and Lemma 2.3 for the 
estimate of solutions in $C^{1+\frac{\gamma}{2}, 2+\gamma}
([0,\www{T}] \times \ooo{\OOO})$.

Extracting a subsequence if necessary, we can assume that $r_m$ converges 
weakly in $H^\gamma(\Omega)$ to $r$ as $m \to +\infty$.
By the Rellich- Kondrashov theorem (see Theorem 1.4.3.2
(p.26) in Grisvard \cite{G}), 
$$
\lim_{m\to+\infty} \Vert r_m - r\Vert_{L^2(\Omega)} = 0.        \eqno{(2.28)}
$$
Furthermore we claim that
$$
\Vert r\Vert_{L^2(\Omega)}\ge \frac{1}{\widetilde C}.      \eqno{(2.29)}
$$
Indeed, by the definition of $r_m$, by means of (2.18), we can obtain 
$$
\Vert r_m\Vert_{H^{\gamma}(\OOO)}=1
\le \widetilde C\left(\Vert r_m\Vert_{L^2(\OOO)} 
+ \frac{\Vert z_{p_m,q_m,a_m,b_m}(T,\cdot)\Vert_{H^{2+\gamma}(\OOO)}
+ \Vert z_{p_m,q_m,a_m,b_m}\Vert_{H^1(\Sigma)} }
{\Vert p_m-q_m\Vert_{H^{\gamma}(\Omega)} } \right).        \eqno{(2.30)}      
$$
Passing to the limit in this inequality as $m\rightarrow +\infty$,
by (2.22) and (2.28) we obtain 
$$
\widetilde C\left(\Vert r_m\Vert_{L^2(\OOO)} 
+ \frac{\Vert z_{p_m,q_m,a_m,b_m}(T,\cdot)\Vert_{H^{2+\gamma}(\OOO)}
+ \Vert z_{p_m,q_m,a_m,b_m}\Vert_{H^1(\Sigma)}}
{\Vert p_m-q_m\Vert_{H^{\gamma}(\Omega)}}\right)
\rightarrow \widetilde C\Vert r\Vert_{L^2(\Omega)}.
$$ 
Hence, (2.30) yields $\www{C}\Vert r\Vert_{L^2(\OOO)} \ge 1$.
Thus (2.29) is proved. 

Passing to the limit in (2.21), by (2.27) we obtain 
$$
\left\{ \begin{array}{rl}
& \ppp_tu_{q,b} = Au_{q,b} - qu_{q,b}, 
                   \quad (t,x) \in (0,+\infty)\times \Omega, \\
& \partial_{\nu_A}u_{q,b}=0 \quad \mbox{on $(0,+\infty)\times \ppp\OOO$},\\
& u_{q,b}(0,\cdot) = b \quad \mbox{in $\OOO$}.
\end{array}\right.  
                                                   \eqno{(2.31)}
$$
We choose $\www{T} = 2T$.  Using the trace theorem and (2.24), we have
$$
\Vert w_m\Vert_{L^2((3\ep,2T)\times \ppp\OOO)} \le C_{\ep} \quad
\mbox{and} \quad
\Vert \ppp_{\nu_A}w_m\Vert_{L^2((3\ep,2T)\times \ppp\OOO)} \le C_{\ep},
\quad m\in \N,
$$
and so we can similarly verify that 
$w_m \longrightarrow w$, $\ppp_{\nu_A} w_m \longrightarrow \ppp_{\nu_A}w$
weakly in $L^2((3\ep,2T)\times \ppp\OOO)$ after extracting subsequences.
By (2.22), we see that $w_m \longrightarrow 0$ in 
$H^1(\Sigma)$ as $m\to +\infty$.  Therefore, we obtain
$w=0$ on $(3\ep, T) \times \Gamma$.  Furthermore 
$\ppp_{\nu_A}w = 0$ on $(3\ep, 2T) \times \ppp\OOO$.
Similarly, in terms of (2.22), we can see that $w(T,\cdot) = 0$ in $\OOO$.

Consequently, we can pass to the limit in (2.23) and by means of
(2.26) and (2.27) we obtain
$$
\left\{ \begin{array}{rl} 
& \ppp_t w = A w - p(x) w + r(x)u_{q,b}(t,x) \quad \mbox{in $(3\ep, 2T)
\times \OOO$},\,\, w(T,\cdot)=0, \\
& \partial_{\nu_A} w = 0 \quad \mbox{on $(3\ep,2T) \times \ppp\OOO$}, 
\quad w = 0 \quad \mbox{on $(3\ep,T) \times \Gamma$}. \\
\end{array}\right.
                                           \eqno{(2.32)}
$$

It is known (see e.g., Pazy \cite{Pa}, Tanabe \cite{Ta}) 
that the function $u_{q,b}(t,\cdot): t\rightarrow L^2(\Omega)$ is analytic 
in $t$ for all $t>0.$  Let $\widetilde \Omega$ be a bounded domain in 
$\Bbb R^n$ with smooth boundary  such that $\Omega\subset \widetilde \Omega$,
$\partial\Omega\setminus\Gamma\subset \partial\widetilde \Omega$ and 
$\widetilde \Omega\setminus\ooo\Omega$ is an open set.  

We make the zero extensions of $w(t,\cdot)$, $u_{q,b}(t,\cdot)$ and $r$ 
to $\widetilde \Omega\setminus\ooo\Omega$. 
We note that after the extensions, 
$w(t,\cdot) \in L^2(\www{\OOO})$ for $3\ep<t<2T$ and
$u_{q,b}(t,\cdot) \in L^2(\www{\OOO})$ for $t>0$.
Moreover, the function $u_{q,b}(t,\cdot): t\rightarrow L^2(\widetilde \Omega)$ 
is analytic in $t>0$.
Hence, the function $ru_{q,b}(t,\cdot): t\rightarrow L^2(\www\Omega)$ 
is analytic in $t>0$. 
Next we extend 
%$r$ in $\widetilde \Omega\setminus\ooo\Omega$ up to 
%a $C^\gamma(\ooo\Omega)$-function 
the coefficients $a_{kj}, b_k,p$ 
on $\widetilde \Omega$ keeping the regularity (1.4) and the positivity 
(1.5).

Let $\chi_{\OOO}(x) 
= \left\{ \begin{array}{rl}
&1, \quad x\in \OOO,\\
&0, \quad x\not\in \OOO.
\end{array}\right.$

Consider an initial boundary value problem
$$
\left\{ \begin{array}{rl} 
& P \www{w}:=\ppp_t \widetilde w -A \widetilde w + p(x)\widetilde  w 
= r(x)\chi_{\Omega}(x)u_{q,b}(t,x) \quad \mbox{in $(3\ep, 2T)\times 
\widetilde \OOO$},\,\,  \\
& \partial_{\nu_A} \widetilde w = 0 \quad \mbox{on $(3\ep,2T) 
\times \ppp\widetilde\OOO$}, \quad\widetilde w(3\ep,\cdot)=w(3\ep,\cdot)
\quad \mbox{in $\www\OOO$}, 
\end{array}\right.
                                           \eqno{(2.33)}
$$ where $w(\ep,\cdot)$ is extended by zero in $\www\Omega\setminus
\ooo\Omega$.
We define 
$\mbox{\bf w}=w $ if $(t,x)\in (3\ep, T)\times \OOO$ and $\mbox{\bf  w=}0$ 
if $(t,x)\in (3\ep, T)\times(\www\Omega\setminus \ooo\OOO)$.
Since $w = \ppp_{\nu_A}w = 0$ on $(3\ep,T) \times \Gamma$ by (2.32), 
it follows that ${\bf w} \in H^{1,2}((3\ep,T) \times \www\OOO)$.
Hence $P\mbox{\bf w}\in L^2((3\ep,T)\times\www\Omega).$ 
On the other hand $P\mbox{\bf w}\vert_{(3\ep,T)\times\Omega}= r u_{q,b}$ and  
$P\mbox{\bf w}\vert_{(3\ep,T)\times(\widetilde\Omega\setminus \Omega)}= 0.$ 
Therefore 
$$
P\mbox{\bf w}= r\chi_{\Omega}u_{q,b}\quad \mbox{on}\,\, 
(3\ep,T)\times\www\Omega
$$
and 
$$
\partial_{\nu_A}\mbox{\bf w}=0\quad \mbox{on}\quad (3\ep,T)\times 
\partial\widetilde \Omega, \quad \mbox{\bf w}(3\ep,\cdot)=w(3\ep,\cdot)
\quad \mbox{in $\www\OOO$}. 
$$
Then the uniqueness of solution to an initial boundary value problem yields
$$
\mbox{\bf w}=\widetilde w\quad \mbox{on}\,\, (3\ep,T)\times\widetilde\Omega.
$$
It is shown in \cite{Kato} and \cite{Kato0} that 
$\widetilde w( t,\cdot): t\rightarrow
 L^2(\widetilde\Omega)$ is analytic in $t$ for all $t>3\ep.$ 
Hence, since ${\bf w} = 0$ in $(3\ep,T) \times (\www\OOO
\setminus \ooo{\OOO})$, we obtain $\widetilde w( t,x)=0$
for all $t> 3\ep$ and $x\in \widetilde \Omega\setminus\Omega.$
Since $\widetilde w\in H^{1,2}((3\ep,2T)\times \Omega)$, we have 
$$
\widetilde w\vert_{(3\ep,2T)\times \Gamma_0}
= \partial_{\nu_A}\widetilde w\vert_{(3\ep,2T)\times \partial\Omega}=0,
$$
where $\Gamma_0=\partial(\widetilde \Omega\setminus\ooo\Omega)\cap
\partial\Omega.$ 
Again the uniqueness of solution to the initial boundary value problem 
implies $w=\widetilde w$ in $(3\ep,2T) \times \OOO$, and we have
$$
w\vert_{(3\ep, 2T)\times \Gamma_0}=0.    \eqno{(2.34)}
$$
Similarly to (2.24), we can verify that $w, \ppp_tw 
\in H^{1,2}((4\ep,2T)\times \OOO)$. 

By Lemma 2.4, there exists a constant $c_0>0$ such that 
$$
u_{q,b}(T,x)\ge c_0 \quad \mbox{for all $x\in \Omega$}.
                                      \eqno{(2.35)}
$$

Then, in terms of (2.34) and (2.35), we apply the uniqueness result 
for the inverse source problem with spatial data at not final time to 
(2.32) (e.g., Theorem 3.1 in Imanuvilov and Yamamoto \cite{IY1}), 
so that we reach 
$r\equiv 0$ in $\OOO$.  This contradicts $r\not\equiv 0$ in (2.29).
Thus (2.22) is false, and the proof of (2.19) is complete.
Thus the proof of Theorem 1 is complete.
$\blacksquare$

\section{Proof of Theorem 2}

We rely on what is called the Reznitskaya transform (e.g., Romanov \cite{Ro}), 
which is the 
main idea for a similar inverse problem solved by Klibanov \cite{K}.

Since function $t\rightarrow u(p)(t,\cdot)\vert_{\partial\Omega}$ and 
$t\rightarrow u(q)(t,\cdot)\vert_{\partial\Omega}$ are analytic in $t>0$ 
(see e.g., Pazy \cite{Pa}, Tanabe \cite{Ta})  in view of 
$u(p) = u(q)$ on $(0,T) \times \Gamma_0$, we obtain
$$
u(p) = u(q) \quad \mbox{on $(0, \infty) \times \Gamma_0$}.
                                       \eqno{(3.1)}
$$
Consider an initial-boundary value problem for a hyperbolic equation:
$$
\ppp_t^2w(t,x) - \Delta w + p(x)w = 0 \quad \mbox{in $(0,\infty)
\times \OOO$},                    \eqno{(3.2)}
$$
$$
\ppp_{\nu_A}w = 0 \quad \mbox{on $(0,\infty) \times \ppp\OOO$}
                                             \eqno{(3.3)}
$$
and
$$
w(\cdot,0) = u_0, \quad \ppp_tw(\cdot,0) = 0 \quad \mbox{in $\OOO$}.
                                                \eqno{(3.4)}
$$
By $u_0 \in H^3(\OOO)$ and $\ppp_{\nu_A}u_0 = 0$ on $\ppp\OOO$,
and $p,q \in W^{1,\infty}(\OOO)$, we can verify 
$$
w(p), \, w(q) \in C^k([0,\infty);H^{3-k}(\OOO)), \quad 
k=0,1,2,3.                                              \eqno{(3.5)}
$$
We can prove (3.5) by Theorem 8.2 (p.275) in \cite{LM}.
Indeed in that theorem, setting $V := H^1(\OOO) 
= \mathcal{D}((-\Delta_N)^{\frac{1}{2}})$ where
$\Delta_Nv = \sum_{j=1}^n \ppp_j^2v$ with 
$\mathcal{D}(\Delta_N) = \{ v\in H^2(\OOO);\, \ppp_{\nu}v\vert
_{\ppp\OOO} = 0\}$, we can accomplish the proof of (3.5) by the theorem and  
the elliptic regularity.
  
Moreover we apply a standard energy estimate and, in view of (1.12), we have 
$$
\Vert w(p)(\cdot,s)\Vert_{L^2(\OOO)} \le C_1\Vert \nabla w(p)\Vert_{L^2(\OOO)},
$$
and we can find a constant $C_2>0$ such that 
$$
\Vert w(p)(\cdot,s)\Vert_{H^1(\OOO)} \le C_2e^{C_2s}\quad 
\mbox{for all $s>0$},
$$
where $C_2>0$ depends on $p, u_0$.  Therefore there exists
$$
\www{u}(p)(t,x):= \int^{\infty}_0 \frac{1}{\sqrt{\pi t}}
e^{-\frac{s^2}{4t}}w(p)(s,x) ds, \quad t>0, \, x\in \OOO,
$$
which is called the Reznitskaya transform in the context of 
inverse problems.

We can prove (e.g., Section 3 of Chapter 6 of Romanov \cite{Ro}) that 
$\www{u}(p)$ satisfies  first two equations  of (1.11) in 
$(0, \infty) \times \OOO$. Let us show that the initial condition also 
holds true.
we remind that $w(p)\in C(\Bbb R^1;L^2(\Omega))$ .
$$
\Vert  \www{u}(p)(t,\cdot)-u_0(0,\cdot)\Vert_{L^2(\Omega)}
= \left\Vert \int_0^{\infty}\frac{1}{\root\of{\pi t}} e^{-\frac{\tau^2}{4t}}
(w(p)(\tau,\cdot)-u_0(\cdot))d\tau\right\Vert_{L^2(\Omega)}
$$
$$
\le \int_0^{\infty}\frac{1}{\root\of{\pi t}} 
e^{-\frac{\tau^2}{4t}}\Vert w(p)(\tau,\cdot)-u_0\Vert_{L^2(\Omega)}d\tau.
$$
For any $\ep>0$, there exists a constant $\delta(\ep)>0$ such that
$$\Vert w(p)(\tau,\cdot)-u_0\Vert_{L^2(\Omega)}
\le \frac{\ep}{2} \quad \mbox{for all $0<\tau < \delta(\ep)$}.
$$
Therefore 
\begin{align*}
&\int_0^{\infty}\frac{1}{\sqrt{t}} e^{-\frac{\tau^2}{4t}}
\Vert w(p)(\tau,\cdot)-u_0\Vert_{L^2(\Omega)}d\tau\\
= &\int_{\delta(\ep)}^{\infty}\frac{1}{\sqrt{t}} 
e^{-\frac{\tau^2}{4t}}\Vert w(p)(\tau,\cdot)-u_0\Vert_{L^2(\Omega)}d\tau
+ \int_0^{\delta(\ep)}\frac{1}{\root\of{\pi t}} e^{-\frac{\tau^2}{4t}}\Vert w(p)(\tau,\cdot)-u_0\Vert_{L^2(\Omega)}d\tau.
\end{align*}
We can estimate 
\begin{align*}
& \int_{\delta(\ep)}^{\infty}\frac{1}{\sqrt{t}} 
e^{-\frac{\tau^2}{4t}}\Vert w(p)(\tau,\cdot)-u_0\Vert_{L^2(\Omega)}d\tau
\le C_2\int_{\delta(\ep)}^{\infty}\frac{1}{\sqrt{t}}
e^{-\frac{\tau^2}{4t}}e^{C_2\tau} d\tau\\
= &C_2 \int_{\delta(\ep)}^{\infty} 2e^{C_2^2t}e^{-\frac{1}{4t}(\tau-2tC_2)^2}
\frac{d\tau}{2\sqrt{t}}
\le 2C_2 e^{C_2^2t} \int_{\frac{1}{2\sqrt{t}}(\delta(\ep)-2tC_2)}^{\infty} 
e^{-\eta^2} d\eta = o(1) \quad \mbox{as $t\to 0$}.
\end{align*}
Consequently we can choose a constants $t_0(\ep)>0$ sufficiently small, so that
$$
\Vert  \www{u}(p)(t,\cdot)-u_0\Vert_{L^2(\Omega)}\le \ep
\quad \mbox{if $0 < t < t_0(\ep)$}.
$$

Thanks to the uniqueness of solution of (1.11), we obtain
$$
u(p)(t,x) = \int^{\infty}_0 \frac{1}{\sqrt{\pi t}}
e^{-\frac{s^2}{4t}}w(p)(s,x) ds, \quad 
u(q)(t,x) = \int^{\infty}_0 \frac{1}{\sqrt{\pi t}}
e^{-\frac{s^2}{4t}}w(q)(s,x) ds, \quad t>0, \, x\in \OOO.
$$
By (3.1), we reach 
$$
\int^{\infty}_0 \frac{1}{\sqrt{\pi t}}e^{-\frac{s^2}{4t}}(w(p)(s,x)-w(q)(s,x)) 
ds=0 \quad t>0, \, x\in \Gamma_0.
$$
Making the change of variables $\widetilde s=e^{-\frac{s^2}{4}}$,
we have 
$$
\int_0^1\tilde s^{\frac 1t-1} \frac{(w(p)-w(q))
(2\,\root\of{-\log \widetilde s},x)}{ \root\of{-\log \widetilde s}}
d\widetilde s=0\quad \mbox{for all $t>0$}.
$$
Setting in the above formula $t=\frac 1n$, choosing $t>0$ such that 
$\frac{1}{t} - 1$ can cover $\N \cup \{0\}$, we obtain 
that the function $\frac{(w(p)-w(q))(2\,\root\of{-\log \widetilde s},x)}
{ \root\of{-\log \widetilde s}}$ is orthogonal to all the polynomials 
on the interval $[0,1]$.
Therefore it  is identically equal to zero by the Weierstrass' polynomial 
approximation theorem.  This implies 
$$
w(p)(\sqrt{\eta},x)-w(q)(\sqrt{\eta},x) = 0 \quad \mbox{for all 
$\eta > 0$ and $x \in \Gamma_0$.}
$$
Hence,
$$
w(p)(t,x) =  w(q)(t,x) \quad \mbox{for all $t>0$ and $x\in \Gamma_0$}.
                                                         \eqno{(3.6)}
$$

Now we reduced the inverse parabolic problem to an inverse hyperbolic problem
for (3.2) - (3.4) with (3.6).
Setting $y:= w(p) - w(q)$, $f:= q-p$ and $R(t,x):= w(q)(t,x)$, we have
$$
\left\{ \begin{array}{rl}
& \ppp_t^2y = \Delta y - p(x)y + f(x)R(t,x), \quad t>0, \, x\in \OOO, \\
& \ppp_{\nu_A}y = 0 \quad \mbox{on $(0,\infty) \times \ppp\OOO$}, \\
& y(0,\cdot) = \ppp_ty(0,\cdot) = 0 \quad \mbox{in $\OOO$}
\end{array}\right.
                                             \eqno{(3.7)}
$$
and
$$
y = 0 \quad \mbox{on $(0,\infty) \times \Gamma_0$.}  
                                                      \eqno{(3.8)}
$$
For direct application of Imanuvilov and Yamamoto \cite{IY2}, we 
will make an extension of $y$ to a wider spatial domain.
First we choose an open smooth domain $\omega \subset \R^n \setminus 
\ooo{\OOO}$ such that $\ooo{\omega} \cap \ppp\OOO \subset \Gamma_0$. 
We fix $T>0$ satisfying $T > \sup_{x\in \OOO} \vert x-x_0\vert$.
Over $\Gamma_0\subset \ppp\OOO$, we take the zero extensions of both
$y$ and $f$ to $\www{\OOO}:= \OOO \cup \omega \cup (\ppp\omega \cap \Gamma_0)$.
By $\ppp_{\nu_A}y = y = 0$ on $(0,T)\times \Gamma_0$, we can readily see that 
$$
y \in C^1([0,T];H^2(\www{\OOO})) \cap C^2([0,T];H^1(\www{\OOO})) 
\cap C^3([0,T];L^2(\www{\OOO})).                      \eqno{(3.9)}
$$
We can extend also $R$ to $\www{\OOO}$ such that 
$$
\vert R(0,x)\vert=\vert u_0(x)\vert\ge \delta > 0 \quad \mbox{on $\ooo{\www{\OOO}}$}, \quad 
R\in H^1(0,T;L^{\infty}(\www{\OOO})).                             \eqno{(3.10)}
$$
Therefore, we can obtain
$$
\left\{ \begin{array}{rl}
& \ppp_t^2y = \Delta y - p(x)y + f(x)R(t,x), \quad t>0, \, x\in \www{\OOO}, \\
& \ppp_{\nu_A}y = 0 \quad \mbox{on $(0,\infty) \times \ppp\www{\OOO}$}, \\
& y(0,\cdot) = \ppp_ty(0,\cdot) = 0 \quad \mbox{in $\www{\OOO}$},\\
& y=0 \quad \mbox{in $(0,T) \times \omega$}. 
\end{array}\right.
$$
In view of the construction of $\omega$ and (1.13), we note that 
$$
\{ x\in \ppp{\www{\OOO}};\, (x-x_0)\cdot \nu(x) \ge 0\} \subset 
\ppp\omega.
$$
Hence, in terms of (3.9) and (3.10),
we can apply the uniqueness for an inverse hyperbolic problem
(e.g., Corollary 3.1 in \cite{IY2}) to conclude that $f=0$ in $\www{\OOO}$,
that is, $p=q$ in $\OOO$.  Thus the proof of Theorem 2 is complete. $\blacksquare$
\section{Proof of Theorem 3}

We define an operator by 
$$
(Kv)(t):= \mu(0) v(t) + \int^t_0 \mu'(t-s)v(s) ds, \quad 0<t<T.
                                                                 \eqno{(4.1)}
$$
Let $y$ satisfy (1.14) and $y(0,\cdot) = 0$ in $\OOO$ and
$y = \vert \nabla y\vert = 0$ on $(0,T)\times \Gamma$.  
Consider an equation with respect to $z(t,x)$:
$$
\ppp_ty(t,x) = (Kz)(t,x), \quad 0<t<T, \, x\in \OOO.
                                                        \eqno{(4.2)}
$$
Since $\mu(0) \ne 0$, the operator $K$ is a Volterra operator of the second 
kind, and we see that $K^{-1}: H^1(0,T) \longrightarrow H^1(0,T)$ exists and 
is bounded.  Therefore, $z(\cdot,x) \in L^2(0,T)$ is well defined for each 
$x \in \OOO$, and
$$
\ppp_ty(t,x) = \mu(0)z(t,x) + \int^t_0 \mu'(t-s)z(s,x) ds, \quad 0<t<T,\, x\in 
\OOO                                                       \eqno{(4.3)}
$$
by $\ppp_ty \in H^{1,2}(Q)$.

Since $y(0,\cdot) = 0$ in $\OOO$ and 
$$
\mu(0)z(t,x) + \int^t_0 \mu'(t-s)z(s,x) ds
= \ppp_t\left( \int^t_0 \mu(t-s)z(s,x) ds\right), 
$$
we obtain
$$
y(t,x) = \int^t_0 \mu(t-s)z(s,x) ds, \quad 0<t<T,\, x\in \OOO,
$$
that is,
$$
y(t,x) = \int^t_0 \mu(s)z(t-s,x) ds, \quad 0<t<T,\, x\in \OOO.
                                                                \eqno{(4.4)}
$$
We will prove that $z \in H^{1,2}(Q)$ satisfies 
$$
\ppp_tz(t,x) = Az(t,x) - p(x)z \quad \mbox{in $(0, \, t_*) \times \OOO$},
                                                                \eqno{(4.5)}
$$
where $t_*\in (0,T)$ is some constant, and 
$$
z = \vert \nabla z\vert = 0 \quad \mbox{on $(0,T) \times \Gamma$}.
                                                           \eqno{(4.6)}
$$
We can readily verify (4.6), because $\ppp_ty = \vert \nabla\ppp_ty\vert 
= 0$ on $(0,T) \times \Gamma$ implies 
$Kz = \vert K(\nabla z)\vert = 0$ on $(0,T) \times \Gamma$, so that  
the injectivity of $K$ directly yields (4.6).

Using $\ppp_ty \in H^{1,2}(Q) \subset C([0,T];L^2(\OOO))$ by (4.3), we have 
$$
\ppp_ty(0,x) = \mu(0)z(0,x), \quad x\in \OOO.
$$
On the other hand, substituting $t=0$ in (1.14) and applying (1.15), we obtain
$$
\ppp_ty(0,x) = \mu(0)f(x), \quad x\in \OOO.
$$
Hence $\mu(0)z(0,x) = \mu(0)f(x)$ for $x\in \OOO$.  By $\mu(0) \ne 0$, we 
reach 
$$
z(0,x) = f(x), \quad x\in \OOO.                  \eqno{(4.7)}
$$

Now we will prove (4.5).  In terms of (4.4) and (4.7), we have
\begin{align*}
& \ppp_ty(t,x) = \mu(t)z(0,x) + \int^t_0 \mu(s)\ppp_tz(t-s,x) ds\\
=& \mu(t)f(x) + \int^t_0 \mu(s)\ppp_tz(t-s,x) ds
\end{align*}
and
$$
(A-p(x))y(t,x) = \int^t_0 \mu(s)(A-p(x))z(t-s,x) ds, \quad 0<t<T,\, x\in \OOO.
$$
Consequently (1.14) implies
\begin{align*}
& \mu(t)f(x) = (\ppp_ty - (Ay- p(x)y)(t,x)\\
=& \mu(t)f(x) + \int^t_0 \mu(s)(\ppp_tz - (A-p(x))z)(t-s,x) ds,
\end{align*}
that is,
$$
\int^t_0 \mu(s)(\ppp_tz - (A-p(x))z)(t-s,x) ds = 0, \quad 0<t<T,\, x\in \OOO.
$$
Hence, setting $Z(s):= \Vert (\ppp_sz - (A-p(x))z)(s,\cdot)\Vert_{L^2(\OOO)}$ 
for $0<t<T$, we reach 
$$
\int^t_0 \mu(s) Z(t-s) ds = 0, \quad 0<t<T.
$$
By the Titchmarsh convolution  theorem (e.g., Titchmarsh \cite{Tit}), 
there exists $t_* \in [0,T]$ such that 
$$
\mu(s) = 0 \quad \mbox{for $0<s<T-t_*$}, \quad
Z(s) = 0 \quad \mbox{for $0<s<t_*$}.
$$
Since $\mu\not \equiv 0$ in $[0,T]$, we see that $T-t_* < T$, that is,
$t_*>0$.  Thus the verification of (4.5) is complete.
$\blacksquare$
\\

Applying the classical unique continuation by (4.5) and (4.6), 
we obtain $z=0$ in $(0,t_*) \times 
\OOO$.  As for the unique continuation, we can refer for example to 
Mizohata \cite{Mi}, Saut-Scheurer \cite{SS} and see also Yamamoto \cite{Ya2009},\cite{Ya2022}.
Thus (4.7) implies $f=0$ in $\OOO.$ The proof of Theorem 3 is complete.
$\blacksquare$

\section{Concluding remarks}

{\bf 5.1.}\\
In this article, we consider two types of inverse problems for parabolic 
equations:
\begin{itemize}
\item
Inverse coefficient problem of determining a 
spatially varying $p(x)$ in 

$$
\left\{ \begin{array}{rl}
& \ppp_tu(t,x) = Au(t,x) - p(x)u(t,x), \quad 0<t<T, \, x\in \OOO,\\
& \ppp_{\nu_A}z = 0 \quad \mbox{on $(0,T)\times \ppp\OOO$},\\
& u(t_0,\cdot) = u_0: \,\, \mbox{given in $\OOO$},
\end{array}\right.
$$
by $u\vert_{(0,T)\times \Gamma}$ with subboundary $\Gamma \subset \ppp\OOO$.
Here $t_0 \in (0,T)$ is given.
\item
Inverse source problem of determining a spatial factor $f(x)$ of 
source terms in
$$
\left\{ \begin{array}{rl}
& \ppp_ty(t,x) = Ay(t,x) - p(x)y(t,x) + \mu(t)f(x), 
\quad 0<t<T, \, x\in \OOO,\\
& y(t_0,\cdot) = u_0: \,\, \mbox{given in $\OOO$},
\end{array}\right.
$$
by $y, \nabla y$ on $(0,T) \times \Gamma$.
\end{itemize}

We summarize our results and the existing results.

{\bf Inverse coefficient problems:}
\begin{enumerate}
\item
Case $0<t_0<T$: there have been many works and among them we can refer for 
example to Beilina and Klibanov \cite{BeKl}, 
Bukhgeim and Klibanov \cite{BK}, Klibanov \cite{K},
Imanuvilov and Yamamoto \cite{IY1} and the references therein.
We refer to Yamamoto \cite{Ya2009, Ya2022} for the stability.  
The main method is based on Carleman estimates.
\item
Case $t_0=T$: The current article solved and established the global 
Lipschitz stability.   We should remark that the uniqueness and also
H\"older stability can be proved by means of the time analyticity
of solution $u$.
Our method can be extended in the case where the coefficients
of the parabolic operator depends on time analytically.
\item
Case $t_0=0$: By  an integral transform,
we can obtain the uniqueness,  
once we have uniqueness for the corresponding hyperbolic 
inverse coefficient problems.  This idea can be found in the book by 
Romanov \cite{Ro}, Klibanov \cite{K}.
By this method, it is impossible to choose any small subboundary where
additional 
data are taken, which seems to be  unnatural for the inverses parabolic problem.
One could prove conditional stability, but it is expected to be of rather weak
rate, which is subject to the transformation method, 
and we do not know whether such 
stability is reasonable by the nature of the inverse problem.
\end{enumerate}

{\bf Inverse source problems:}
The main concern is the uniqueness and we should expect 
the uniqueness without boundary condition on the whole
lateral boundary $(0,T) \times \ppp\OOO$.
\begin{enumerate}
\item
Case $0<t_0<T$: The method by Carleman estimate gives the uniqueness and also 
conditional stability in some subdomain $\subset \OOO$ if $\mu(t_0) \ne 0$.
\item
Case $t_0=T$: still open.
\item
Case $t_0=0$: The current article establishes the uniqueness if 
$\mu(0) \ne0$.  However under more generous but reasonable condition
$\mu \not\equiv 0$ on $[0,T]$, we do not know the uniqueness in general.
\end{enumerate}

\vspace{0.2cm}

{\bf 5.2.}
For the case $t_0=0$, our available method seems only the Reznitskaya 
transform.  Within adequate regularities, the transform transfers 
the uniqueness for inverse hyperbolic problems to the uniqueness for 
inverse parabolic problem.
Thus the primary technical highlight is inverse hyperbolic problems.

When one can obtain the uniqueness 
for more general inverse hyperbolic problems, we can make it 
produce the uniqueness for the corresponding inverse parabolic problems
not only in the case of $A = \Delta$.

In order to treat general $A$, we need some restrictive assumptions for 
the principal coefficients $a_{jk}$ of $A$.
Now we formulate such assumptions as follows.
We set $i:= \sqrt{-1}$, $Q:= (0,T)\times \OOO$, 
$\nabla:= (\ppp_1,..., \ppp_n)$,
$\nabla_{t,x}:= (\ppp_t, \ppp_1, ..., \ppp_n)$,
$\nabla_{\xi}:= (\frac{\ppp}{\ppp \xi_0}, \frac{\ppp}{\ppp \xi_1}, ...,
\frac{\ppp}{\ppp \xi_n})$ for $\xi:= (\xi_0, ..., \xi_n)$. For any $C^1$ functions $g(x,\xi) $ and $f(x,\xi)$ we introduce the Poisson bracket of these functions 
$$
\{f,g\}=\sum_{j=1}^n \partial_{\xi_j}f\partial_{x_j}g-\partial_{x_j}f\partial_{\xi_j} g+\partial_{\xi_0}f\partial_{t}g-\partial_{t}f\partial_{\xi_0}g.
$$

For the  hyperbolic operator:
$$
\widetilde P(x,D)v: =\partial^2_tv
- \sum_{j,k=1}^n \partial_j(a_{jk}(x)\partial_kv)
+ \sum_{j=1}^n b_j(x)\partial_jv + p(x)v
$$
we define principal symbol of this operator  $\www{p}(x,\xi)$ by  formula
$$
\www{p}(x,\xi):=\xi_0^2-a(x,\xi,\xi),\quad x\in \OOO,\,
\xi=(\xi_0,\dots,\xi_n).
$$

First we 
assume that a function $\www \psi(t,x)$ is pseudo-convex with respect to 
$\widetilde P$ in $Q$, that is,
$$
\left\{ \begin{array}{rl}
& \nabla_{t,x}\www\psi(t,x)\ne 0\quad \mbox{for all $(t,x)\in \overline{Q}$},
                              \cr\\
& \frac{ \{ \www{p}(x,\xi-i\tau\nabla_{t,x}\widetilde \psi(t,x)),\,\,
\www{p}(x,\xi+i\tau \nabla_{t,x}\widetilde \psi(t,x)) \}}{2i\tau}>0 \cr\\
& \mbox{for all $(t,x,\xi)\in Q\times(\Bbb R^{n+1}\setminus\{0\})$
and $\tau>0$ satisfying}\cr\\
& \www{p}(x,\xi+i\tau \nabla_{t,x}\www\psi(t,x))
\, =0.
\end{array}\right.
                                          \eqno{(5.4)}
$$
Second we assume that a subboundary $\Gamma_1 \subset \ppp\OOO$ satisfies
$$
a(x,\nu(x), \nabla\www \psi(t,x)) < 0 \quad \mbox{for all
$(t,x)\in [0, T]\times \overline{\partial\Omega\setminus\Gamma_1}$ },
                                      \eqno{(5.5)}
$$
and
$$
\partial_t\www\psi(T,x)<0 \quad \mbox{and}\quad \partial_t\www\psi(0,x)>0\quad 
\mbox{for all $x\in \ooo\Omega$}.                \eqno{(5.6)}
$$
Finally we assume that for each $y\in \overline{\partial\Omega\setminus 
\Gamma_1}$, there exist an open neighborhood $U(y)$ of $y$ and a function 
$\psi_1$ defined in $U(y)$ such that $\psi_1$ is pseudo-convex 
with respect to the symbol $\www{P}$ in $(0,T) \times U(y)$ and 
%$$
%\psi_1(x)\le \psi(x)\quad \mbox{for all $x\in U(y)$}  \eqno{(5.7)}
%$$
%and 
$$
a(x,\nu(x), \nabla\psi_1(x))>0 \quad \mbox{for all 
$x\in \partial\Omega\cap U(y)$}.                    \eqno{(5.7)}
$$

Under the assumptions (5.4) - (5.7), replacing the subboundary $\Gamma_0$ 
defined in (1.13) by $\Gamma_1$ and choosing sufficiently large 
$T>0$ for the principal coefficients $a_{jk}(x)$, 
we can generalize Theorem 2 and obtain the same uniqueness in determining 
the zeroth-order coefficient $p(x)$ for (1.1)-(1.3) with $t_0=0$
by extra data $u\vert_{\Gamma_1\times (0,T)}.$

The conditions (5.4) - (5.7) are satisfied for example in the case
$a_{jk}=\delta_{jk}$  but in general
it is not at all simple to find a generous sufficient condition on 
$a_{jk}$ with $1\le j,k \le n$ admitting $\www\psi$ and $\psi_1$ 
which satisfy (5.4) - (5.7).
For example, Imanuvilov and Yamamoto \cite{IY03},
\cite{IYLame} discuss  such generous conditions respectively for 
the cases of the Dirichlet and the Neumann type of boundary conditions.
However we do not exploit more.
\\

{\bf Acknowledgement}\\ 
M.\! Yamamoto is supported by Grant-in-Aid for Scientific Research (A) 
20H00117 and Grant-in-Aid for Challenging Research (Pioneering) 21K18142, JSPS.

\end{document}